\documentclass[a4, 12pt, 
]{amsart}

\usepackage[all]{xy}
\usepackage{mathrsfs, amsthm, amscd, amsmath, amssymb, graphicx}
\usepackage{float}
\usepackage{graphicx}

\usepackage[dvipdfmx]{hyperref}

\newtheorem{theo}{Theorem}[section]
\newtheorem{prop}[theo]{Proposition}
\newtheorem{lemm}[theo]{Lemma}

\newtheorem{prob}[theo]{Problem}
\newtheorem{conj}[theo]{Conjecture}

\numberwithin{equation}{section}

\theoremstyle{definition}

\theoremstyle{remark}

\newcommand{\Alb}[0]{\operatorname{Alb}}

\newcommand{\Sym}[0]{\operatorname{Sym}}

\newcommand{\deldel}{\sqrt{-1}\partial \overline{\partial}}

\newcommand{\reg}{{\rm{reg}}}

\usepackage{geometry}
\newenvironment{dedication}
  {
   \itshape             
   \centering          
  }

\begin{document}

\title[Open problems on  positively curved projective varieties]
{Open problems on structure \\of positively curved projective varieties}

\author{Shin-ichi MATSUMURA}

\address{Mathematical Institute, Tohoku University, 
6-3, Aramaki Aza-Aoba, Aoba-ku, Sendai 980-8578, Japan.}

\email{{\tt mshinichi-math@tohoku.ac.jp}}
\email{{\tt mshinichi0@gmail.com}}

\date{\today, version 0.01}

\renewcommand{\subjclassname}{%
\textup{2010} Mathematics Subject Classification}


\maketitle

\begin{dedication}
Dedicated to Professor Ahmed Zeriahi on the occasion of his retirement  
\end{dedication}

\begin{abstract}
We provide supplements and open problems related to structure theorems for maximal rationally connected fibrations of certain positively curved projective varieties, 
including smooth projective varieties with semi-positive holomorphic sectional curvature, pseudo-effective tangent bundle, and nef anti-canonical divisor. 
\end{abstract}



\section{Introduction} \label{Sec1}

Certain ``positively curved" varieties, 
which are often formulated 
to have positive holomorphic bisectional curvatures, tangent bundles, or anti-canonical divisors, 
have occupied an important place in the classification theory of projective varieties. 
The Frankel conjecture in differential geometry and the Hartshorne conjecture in algebraic geometry, 
proved by Siu-Yau in \cite{SY80} and Mori in \cite{Mor79}, 
have given a beautiful characterization of projective spaces, 
respectively, 
in terms of positive holomorphic bisectional curvatures and ample tangent bundles.
Since that time, it has become clear that 
structures of positively curved varieties are closely related to the geometry of rational curves 
and more restricted than those of negatively curved varieties (i.e., they have a certain rigidity); 
for example, Fano manifolds are always rationally connected 
(i.e., any two points can be connected by a rational curve),   
and also the Albanese map of compact K\"ahler manifolds with semi-positive holomorphic bisectional curvature is locally trivial (i.e., all the fibers are isomorphic each other). 
One of the central problems in this field 
is to understand structures reflecting rational curves and rigidities, 
by using naturally associated fibrations, such as Albanese maps, Iitaka fibrations, and maximal rationally connected fibrations.

In this paper, we mainly study maximal rationally connected fibrations (MRC fibrations for short) of projective varieties. 
An MRC fibration $\phi: X \dashrightarrow Y$ of a projective variety $X$, 
introduced in \cite{Cam92, KoMM92}, 
is a rational map to a certain negatively curved variety $Y$ 
whose general fiber $F$ is compact and rationally connected 
(see Section \ref{Sec2} for the precise definition and properties). 
The crucial structure theorems have been established for MRC fibrations of a smooth projective variety 
$X$ with semi-positive holomorphic bisectional curvature (resp. nef tangent bundle, nef anti-canonical divisor) 
in  \cite{HSW81, Mok88} (resp. \cite{DPS94}, \cite{Cao19, CH19}); 
for example, under the above positivity assumption, 
the variety $X$ admits a holomorphic and locally trivial MRC fibration 
$\phi:X \to Y$  onto a certain flat manifold $Y$, 
from which $X$ can be decomposed into the rationally connected fiber $F$ and the certain flat base $Y$.
In this paper, we give supplements and open problems related to structure theorems for MRC fibrations, 
reviewing recent  generalizations of the results in \cite{HSW81, Mok88} 
(resp. \cite{DPS94}, \cite{Cao19, CH19}) 
to holomorphic sectional curvatures (resp. pseudo-effective tangent bundles, 
nef anti-log canonical divisors of klt pairs). 

The remainder of this paper is organized as follows: 
In Section \ref{Sec2},  we recall the notions and properties of rationally connected varieties and MRC fibrations. 
Moreover, we explain a technique in the theory of (holomorphic) foliations to take a holomorphic MRC fibration, 
which is one of the key points in establishing structure theorems for MRC fibrations. 
In Section \ref{Sec3}, we consider results for holomorphic sectional curvatures 
in \cite{Mat18a, Mat18b}, 
which generalize the results of \cite{HSW81, Mok88} and Yau's conjecture on positive holomorphic sectional curvature. 
In Section \ref{Sec4}, we focus on studies of pseudo-effective tangent bundles initiated in \cite{HIM19}, 
explaining differences from the structure theorem for nef tangent bundles in \cite{DPS94}. 
In Section \ref{Sec5}, we introduce results for projective klt pairs with nef anti-log canonical divisor in \cite{CCM19}, 
with the goal of generalizing the results in \cite{Cao19, CH19} to log pairs.

\subsection*{Acknowledgements}
The author is grateful to an anonymous referee for his/her helpful comments 
and giving the reference \cite{Par11}.
He is supported by the Grant-in-Aid 
for Young Scientists (A) $\sharp$17H04821 and 
Fostering Joint International Research (A) $\sharp$19KK0342 from JSPS. 

\textit{It is the great pleasure of the author to dedicate this paper to Professor Ahmed Zeriahi, 
in honor of his outstanding achievements and 
with memories of their time together as lecturer at SEAMS School in Hanoi. }

\section{Maximal Rationally Connected Fibrations}\label{Sec2}

In this section, we review the basic notions and their properties related to the geometry of rational curves, 
including uniruled varieties, rationally connected varieties, and MRC fibrations.

A curve is called a \textit{rational curve} if its normalization is the one-dimensional projective space $\mathbb{P}^1$. 
A projective variety $X$ is said to be \textit{uniruled} if $X$ is covered by rationals curves. 
In the case of $X$ being smooth, 
it follows from \cite{BDPP} that $X$ is uniruled if and only if the canonical divisor $K_X$ is not pseudo-effective. 
A projective variety $X$ is said to be \textit{rationally connected} (resp. {\textit{rationally chain connected}})
if any two points can be connected by one rational curve (resp. a chain of rational curves). 
In the case of $X$ having at worst dlt singularities, 
it follows from \cite[Corollary 1.5]{HM07}  that 
the rational connectedness is equivalent to the rational chain connectedness. 

MRC fibrations of projective varieties, introduced in \cite{Cam92, KoMM92},  interpolate the rational connectedness and uniruledness. 
We recall that the definition of MRC fibrations of a projective variety $X$ with mild singularities (e.g., with klt singularities). 
A rational map $\phi:X \dashrightarrow Y$  is called an \textit{RC fibration} (resp. \textit{MRC fibration} )
if it satisfies the first two conditions (resp. all the conditions) below: 
\begin{itemize}
\item $\phi:X \dashrightarrow Y$ is an almost holomorphic map to a projective variety $Y$ (i.e., a rational map whose general fibers are compact). 
\item General fibers of $\phi:X \dashrightarrow Y$ are rationally connected. 
\item There is no horizontal rational curve (i.e., no rational curve whose image under $\phi$ is not one point) 
passing through a general point in $X$. 
\end{itemize}
MRC fibrations are trivial in extreme cases; 
the identity map ${\rm{id}}_X: X \to X$ of $X$ is an MRC fibration if and only if $X$ is not uniruled, 
and also the constant map $X \to \{\rm{1pt}\}$ is an MRC fibration  if and only if $X$ is rationally connected. 
MRC fibrations are not uniquely determined by a given projective variety $X$; 
there is an ambiguity in the choices of the birational models of its image $Y$. 
Thanks to this ambiguity, we may assume that the base $Y$ is smooth by taking a resolution of the singularities of $Y$.

A typical example of RC fibrations is the projective space bundle $\mathbb{P}(E) \to Y$ 
associated with a (holomorphic) vector bundle  $E$ on a smooth projective variety $Y$.  
For example, when $Y$ is an abelian variety, 
the total space $\mathbb{P}(E)$ has no horizontal rational curve, 
and thus $\mathbb{P}(E) \to Y$ is an MRC fibration. 
However, when $Y$ is the projective space $\mathbb{P}^k$, the morphism $\mathbb{P}(E) \to Y$ is not an MRC fibration. 
In general, for an RC fibration $\phi:X \dashrightarrow Y$, 
any rational curves in $Y$ can be lifted into $X$ by \cite{GHS03}, 
and hence the third condition in the above definition can be rephrased as follows: 
the image $Y$ is not uniruled; 
equivalently, the canonical divisor $K_Y$ of $Y$ is pseudo-effective by \cite{BDPP} 
when $Y$ is smooth.

MRC fibrations are almost holomorphic by definition, 
but they are not necessarily represented by holomorphic maps. 
In fact, there is a projective variety $Z$ with Picard number one, 
such that $Z$ is uniruled but not rationally connected, 
which was constructed by Fujita (see \cite[Example 6.5]{EIM20} for details). 
The variety $Z$ admits no non-trivial fibrations 
(in particular, no holomorphic MRC fibrations) 
since the Picard number of $Z$ is one. 
Note that $Z$ has bad singularities: 
$Z$ is neither log canonical nor $\mathbb{Q}$-factorial.
The author does not know an example of smooth projective varieties (or varieties with mild singularities) admitting no holomorphic MRC fibration. 
Now we suggest the following problem:

\begin{prob}
When do we choose an MRC fibration to be holomorphic from the ambiguities in the choices of $Y$?
\end{prob}

This fundamental problem plays a crucial role when we study structure theorems for MRC fibrations. 
The following lemma, based on the Reeb stability, gives a useful sufficient condition to represent MRC fibrations by a holomorphic and smooth morphism. 

\begin{lemm}[{\cite[Corollary 2.11]{Hor07}}]\label{fol-lem}
Let $X$ be a compact K\"ahler manifold and $W \subset T_X$ be an integrable subbundle of $X$. 
If the foliation $W$ has at least one compact and rationally connected leaf,  
then there exists a smooth morphism $X \to Z$ such that 
$W$ coincides with the relative tangent bundle $T_{X/Z}$. 
\end{lemm}

When $X$ is a smooth projective variety, 
we can take a non-empty Zariski open set $Y_0 \subset Y$, so that 
an MRC fibration $\phi: X \dashrightarrow Y$ is a smooth morphism over $Y_0$ 
since $\phi: X \dashrightarrow Y$ is almost holomorphic.  
Then, the relative tangent bundle $T_{X/Y}$ can be defined as a subbundle of $T_X$ on the inverse image $X_{0}:=\phi^{-1}(Y_0)$. 
A general leaf of $T_{X/Y}$ (defined only on $X_0$) is compact and rationally connected by definition. 
Therefore, if this subbundle on $X_0$ can be extended to a subbundle $W \subset T_{X}$ on $X$, 
then the smooth morphism $X \to Z$ obtained from the lemma gives a holomorphic MRC fibration of $X$. 
This strategy actually works in the proof of the structure theorems 
introduced in Section \ref{Sec3}, \ref{Sec4}, and \ref{Sec5}, 
interestingly enough for different reasons.

\section{Semi-Positive Holomorphic Sectional Curvature} \label{Sec3}

The Frankel conjecture proved by Siu-Yau in \cite{SY80} says that 
any smooth projective varieties admitting a K\"ahler metric with positive holomorphic bisectional curvature are isomorphic to the projective space. 
As one of the extensions of the Frankel conjecture, 
Howard-Smyth-Wu and Mok established the structure theorem for compact K\"ahler manifolds $M$ with semi-positive holomorphic bisectional curvature; 
Howard-Smyth-Wu in \cite{HSW81} proved that $M$ admits a locally trivial morphism $f: M\to B$,  
so that  the image $B$ has the flat tangent bundle and 
that the fiber $F$ admits a K\"ahler metric 
whose holomorphic bisectional curvature is semi-positive and whose Ricci curvature is quasi-positive. 
This implies that the fiber $F$ is a Fano manifold (in particular, rationally connected), and thus $f: M\to B$ is automatically an MRC fibration of $M$. 
Moreover, Mok in \cite{Mok88} proved that $F$ satisfying the above conditions is a Hermitian symmetric manifold.

This subsection is devoted to explaining some recent progresses of semi-positive \textit{holomorphic sectional curvature}. 
The holomorphic bisectional curvature ${\rm{BSC}}_g$ and sectional curvature ${\rm{HSC}}_g$ of a K\"ahler metric $g$ 
are defined to be 
$$
{\rm{BSC}}_g(v,w):=\dfrac{R_g(v, \bar v,w, \bar w)}{|v|_g^2 |w|_g^2} \quad \text{ and } \quad 
{\rm{HSC}}_g(v):=\dfrac{R_g(v, \bar v,v, \bar v)}{|v|_g^4 }
$$
for (non-zero) tangent vectors $v, w \in T_X$, 
where $R_g$ is the curvature tensor associated with $g$. 
The holomorphic sectional curvature ${\rm{HSC}}_g$ determines the curvature tensor $R_g$ 
(which means that, if ${\rm{HSC}}_g={\rm{HSC}}_h$ for K\"ahler metrics $g$ and $h$, then we have $R_g=R_h$), 
but there is no explicit relation between ${\rm{HSC}}_g$ and ${\rm{BSC}}_g$. 
Hence it is interesting to pursue an analogy or a difference 
between holomorphic sectional curvature and bisectional curvature. 

The positivity of holomorphic sectional curvatures is much weaker than that of holomorphic bisectional curvature. 
In fact, it follows from Hitchin's result in \cite{Hit75} that 
the Hirzebruch surfaces have positive holomorphic sectional curvature 
(see \cite{AHZ18}  for a generalization of Hitchin's result).
This  tells us that smooth projective varieties with positive holomorphic sectional curvature
is not necessarily Hermitian symmetric and not even Fano. 
This example is in contrast to the case of  semi-positive bisectional curvature 
or negative holomorphic sectional curvature. 
In fact, if a compact K\"ahler manifold $X$ has semi-positive bisectional curvature 
(more generally, the nef tangent bundle), 
then $X$ contains no submanifold $Z$ with negative normal bundle, 
which easily follows from the standard exact sequence: 
$$
0 \to T_{Z} \to T_{X}|_{Z} \to N_{X/Z} \to 0. 
$$
In particular, when $\dim X =2$, the surface $X$ has no negative curves.  
Yau's conjecture on negative holomorphic sectional curvature and its solution in \cite{DT19, TY17, WY16}
asserts that smooth projective varieties with negative holomorphic sectional curvature 
have the ample canonical divisor. 
Nevertheless, the Hirzebruch surfaces 
except for $\mathbb{P}^{1} \times \mathbb{P}^{1}$ 
have a negative curve and are not Fano.

The following conjecture posed by Yau can be regarded as 
an analogy of Mok's result for holomorphic sectional curvature. 

\begin{conj}[Yau's conjecture, {\cite[Problem 47]{Yau82}}]\label{conj-Yau}
If a compact K\"ahler manifold $X$ has positive holomorphic sectional curvature, 
then $X$ is projective and rationally connected.  
\end{conj}

Yang in \cite{Yan18a} solved Yau's conjecture by introducing the notion of the RC positivity 
(see also \cite{Yan18b, Mat13, Yan19} and references therein for the RC positivity). 
In the additional assumption of $X$ being projective, 
Heier-Wong in \cite{HW15} generalized Yau's conjecture to quasi-positive holomorphic sectional curvatures. 
Yau's conjecture was further generalized in \cite{Mat18b}, 
by using the invariant $n_{{\rm{tf}}}{(X, g)}$ defined by 
$$
n_{{\rm{tf}}}{(X, g)}:=\dim X - \inf_{p \in X} \dim V_{{\rm{flat}},p},  
$$
where $V_{{\rm{flat}},p}$ is the subspace of the tangent space $T_{X,p}$ at $p$
consisting of all the truly flat tangent vectors $v$ introduced in \cite{HLWZ18}. 
Here a tangent vector $v \in T_{X,p}$ at a point $p \in X$ is said to be {\textit{truly flat}} if 
it satisfies that 
$$R_g(v, \bar x, y, \bar z)=0
$$ for any $x,y,z \in T_{X,p}$. 
The invariant $n_{{\rm{tf}}}{(X, g)}$ can be seen as an analog of the numerical Kodaira dimension 
(see \cite{Nak04} for the definition) 
and measures the positivity of holomorphic sectional curvatures.

\begin{theo}[{\cite[Theorem 1.2]{Mat18b}}]
\label{thm-hsc-posi} 
Let $X$ be a smooth projective variety and  $g$ be
a K\"ahler metric  with semi-positive holomorphic sectional curvature. 
Let  $\phi: X \dashrightarrow Y$ be an MRC fibration of $X$. 
Then we obtain 
$$
\dim X -\dim Y \geq n_{{\rm{tf}}}{(X, g)}. 
$$
\end{theo}
If the holomorphic sectional curvature is quasi-positive, 
it can be seen that $n_{{\rm{tf}}}{(X, g)}=\dim X$, 
and thus $\dim Y=0$ by the above theorem, 
which means that $X$ is rationally connected. 
This theorem is expected to still hold in the case where $X$ is a compact K\"ahler manifold.
In fact, when a compact K\"ahler manifold $X$ has positive holomorphic sectional curvature, 
Yang proved that $h^{0}(X, \Omega_X^q)=0$ for any $q>0$, and thus $X$ is automatically projective (see \cite[Theorem 1.7]{Yan18a}). 
Then, the following problem naturally arises as a generalization of Yang's criteria for projectivity.

\begin{prob}[{\cite[Problem 5.1]{Mat18b}}]
\label{prob-projective}
Let $(X, g)$ be a compact K\"ahler manifold 
with  semi-positive holomorphic sectional curvature. 
Assume that $n_{{\rm{tf}}}{(X, g)} = \dim X$. 
\begin{itemize}
\item[$\bullet$] Does it hold that $h^{0}(X, \Omega_X^q)=0$ for any $q>0$? 
\item[$\bullet$] Is $X$ automatically projective? 
\end{itemize}
\end{prob}

By combining a technique in the proof of Theorem \ref{thm-hsc-posi} with the theory of foliations, 
we can obtain the following structure theorem for semi-positive holomorphic sectional curvatures.

\begin{theo}[{\cite[Theorem 1.3]{Mat18b}}]
\label{thm-hsc-semi}
Let $X$ be a smooth projective variety and  $g$ be
a K\"ahler metric  with semi-positive holomorphic sectional curvature. 
Then we obtain\,$:$
\begin{itemize}
\item[$\bullet$]
There exists a surjective morphism $\phi: X \to Y$ 
to a smooth projective variety $Y$ with the following properties\,$:$
\begin{itemize}
\item The morphism $\phi: X \to Y$ is a locally trivial morphism. 
\item The image $Y$ is a smooth projective variety with a flat metric. 
In particular, there exists a finite \'etale cover $A \to Y$ by an abelian variety $A$. 
\item The fiber $F$ is a rationally connected manifold. 
In particular, the morphism $\phi: X \to Y$ is an MRC fibration of $X$. 
\end{itemize}
In particular, the fiber product $X^*:=A \times_Y X$ admits 
the locally trivial Albanese map $X^* \to A$ to the abelian variety $A$ 
with the rationally connected fiber $F$. 

\item[$\bullet$]
We obtain an isomorphism 
$$
X_{\rm{univ}} \cong \mathbb{C}^m \times F, 
$$
where $X_{\rm{univ}}$ is the universal cover of $X$, 
and $F$ is the rationally connected fiber of $\phi$. 
Moreover, there exists a representation $\rho : \pi_1(Y) \to {\rm{Aut}}(F)$ such that 
$X$ is isomorphic to $\mathbb{C}^m \times F/\pi_1(Y)$.

\item[$\bullet$] 
In particular, the fundamental group of $X$ is an extension of a finite group by $\mathbb{Z}^{\oplus 2m}$. 

\item[$\bullet$] There exist a K\"ahler metric $g_F$ on the fiber $F$ and 
a K\"ahler metric $g_Y$ on $Y$ with the following properties\,$:$

\begin{itemize}
\item The holomorphic sectional curvature of $g_F$ is semi-positive.  
\item The K\"ahler metric $g_Y$ is flat. 
\item The above isomorphism $X_{\rm{univ}} \cong \mathbb{C}^m \times F$ is not only biholomorphic but also isometric 
with respect to the K\"ahler metrics $\mu^*g $, $\pi^*g_{Y}$, and $g_F$. 
\end{itemize}
Here $\pi$ and $\mu$ respectively denote the universal cover $\pi: \mathbb{C}^m \to Y$ of $Y$ 
and the universal cover $\mu: X_{\rm univ} \to X$ of $X$. 
\end{itemize}
\end{theo}
We have the following commutative diagram$:$
\begin{equation*}
\xymatrix@C=40pt@R=30pt{
 X_{\rm univ} \cong \mathbb{C}^{m} \times F\ar[r]^{} \ar[d]^{} & X^*:=A \times_Y X\ar[d] \ar[r]^{} & X\ar[d]^{\phi} \\ 
\mathbb{C}^{m} \ar[r]  & A \ar[r]^{}  &  Y.\\   
}
\vspace{0.1cm}
\end{equation*}

The morphism $\phi: X \to Y$ in the theorem can be obtained 
from the isomorphism $X_{\rm{univ}} \cong \mathbb{C}^m \times F$ 
and the representation $\rho : \pi_1(Y) \to {\rm{Aut}}(F)$ as mentioned above, 
which is a stronger property than the local triviality of $\phi$. 
This implies that the projective space bundle over an elliptic curve 
does not necessarily admit semi-positive holomorphic sectional curvature, 
whereas, if a smooth projective variety $Y$ has positive holomorphic sectional curvature, 
then any projective space bundles over $Y$ does so by \cite{AHZ18}. 
Indeed, the projective space bundle  associated with  
the vector bundle $\mathcal{O}_Y \oplus \mathcal{O}_Y(np)$ 
is not constructed by the representation of the fundamental group $\pi_{1}(Y) \cong \mathbb{Z}^{\oplus 2}$, 
where $Y$ is an elliptic curve, $p \in Y$, and $n \in \mathbb{Z}$. 
This negatively answers a question posed in \cite{AHZ18}.

It is an attractive problem to generalize the above structure theorem to compact K\"ahler manifolds. 

\begin{prob}\label{prob-m}
Can we generalize Theorem \ref{thm-hsc-semi} to  compact K\"ahler manifolds?
\end{prob}

The fundamental group plays an important role 
when we consider such a structure theorem. 
For a smooth projective variety $X$ with semi-positive holomorphic sectional curvature, 
we can conclude that 
the fundamental group of $X$ is an extension of a finite group 
by a free abelian group by the structure theorem. 
The same conclusion can be expected even when $X$ is a compact K\"ahler manifold, 
but it is still an open problem. 

We now suggest a strategy to solve Problem \ref{prob-m}. 
Let us consider the Albanese map of $X$ (or an \'etale covering space of $X$) instead of MRC fibrations. 
It is not so difficult to check the same  conclusion as in Theorem \ref{thm-hsc-semi} 
for the Albanese map except for the rational connectedness of fibers. 
Then, the fiber $(F, g_F)$ can be expected to satisfy a certain quasi-positivity (e.g.,  $n_{{\rm{tf}}}{(F, g_F)} = \dim F$),  
compared to the structure theorem of Howard-Smyth-Wu and Mok. 
Together with Problem \ref{prob-projective}, 
we can expect that the fiber $F$ is projective and rationally connected.  
We summarize the above problems in the following form:

\begin{prob}[cf. {\cite[Problem 5.2]{Mat18b}}]
\label{prob-main}
Let $(X, g)$ be a compact K\"ahler manifold 
with semi-positive holomorphic sectional curvature. 

\begin{itemize}
\item[$\bullet$] Is the fundamental group of $X$ is an extension of a finite group 
by $\mathbb{Z}^{\oplus 2m}$?

\end{itemize}
For  an appropriate finite \'etale cover $X^* \to X$, 
we consider its Albanese map $\alpha: X^* \to \Alb(X^*)$. 
\begin{itemize}
\item[$\bullet$] Does the fiber $F$ admit a K\"ahler metric $g_F$ such that $n_{{\rm{tf}}}(F, g_F)=\dim F$?
\item[$\bullet$] Is the fiber $F$  projective and rationally connected? 
\end{itemize}
\end{prob}

Let us observe the proof of Theorem \ref{thm-hsc-semi}  in detail. 
It is shown in the proof that the tangent bundle $T_X$ splits into the direct sum: 
$$
T_X \cong T_{X/Y} \oplus \phi^{*} T_Y. 
$$
Moreover, it is also shown that the subbundle $\phi^{*} T_Y \subset T_X$ is integrable and 
all the tangent vectors in  $\phi^{*} T_Y \subset T_X$ are truly flat. 
For compact K\"ahler manifolds with semi-positive holomorphic sectional curvature, 
it is quite interesting to ask whether the converse implication holds, 
more specifically, whether truly flat tangent vectors determine a foliation and its foliation induces an MRC fibration. 
Such a question is also interesting when we consider semi-negative holomorphic sectional curvatures. 
Indeed,  Heier-Lu-Wong-Zheng in \cite{HLWZ18} proved that  the abundance conjecture leads to 
the structure theorem for the Iitaka fibration of smooth projective varieties with semi-negative holomorphic sectional curvature. 
Further, in their proof, a similar splitting theorem of the tangent bundle was obtained. 
Hence, if the same strategy works (i.e., truly flat tangent vectors determine a foliation), 
then its foliation may directly construct the Iitaka fibration and 
solve the abundance conjecture in the case of semi-negative holomorphic sectional curvatures. 

\begin{prob}\label{prob-main}
Let $(X, g)$ be a compact K\"ahler manifold 
with  semi-positive $($or semi-negative$)$ holomorphic sectional curvature. 
Let  $V \subset T_X$ be the set of ``appropriate"  truly flat tangent vectors. 

\begin{itemize}
\item[$\bullet$] Does $V$ determine an integrable subbundle of  $T_X$?
\item[$\bullet$] Can we find a subbundle $W \subset T_X$ such that $T_X \cong W \oplus V $?
\item[$\bullet$] In the case of semi-positive holomorphic sectional curvatures, 
can we construct an MRC fibration $f: X \to Y$ such that $W=T_{X/Y}$?
\item[$\bullet$] In the case of semi-negative holomorphic sectional curvatures, 
can we construct the Iitaka fibration $f: X \to Z$ associated with the canonical bundle $K_X$
such that $V=T_{X/Z}$?
\end{itemize}

\end{prob}

The notion of nef tangent bundles is an analog  in algebraic geometry of semi-positive holomorphic bisectional curvatures 
(see \cite{DPS94} and Section \ref{Sec4} for nef tangent bundles). 
Hence, it is of interest to ask what an analog of semi-positive holomorphic sectional curvature is. 
It is also a fundamental problem to classify all the varieties of fixed dimension admitting (semi-)positive holomorphic sectional curvature, 
by applying Theorem \ref{thm-hsc-semi}. 
The complete classification is not known even in the case of surfaces. 
One of the difficulties is the existence of negative curves. 
When the tangent bundle $T_{X}$ of $X$ is nef, 
there is no negative curve in $X$, 
and hence the blow-ups of varieties can be excluded. 
When we consider the positive holomorphic sectional curvature, 
we can not exclude the case of $X$ having a negative curve; 
indeed, the Hirzebruch surfaces have positive sectional curvature. 
It is not even obvious whether the blow-up of Hirzebruch surfaces has 
a positive holomorphic sectional curvature.

\begin{prob}[cf.{\cite[Problem 5.3]{Mat18b}}]
\label{prob-cl}
\ 
\begin{itemize}
\item[$\bullet$] Does the blow-up of Hirzebruch surfaces admit  
a positive holomorphic sectional curvature?

\item[$\bullet$] Can we classify all the varieties of fixed dimension $($e.g., dimension two or three$)$ 
with semi-positive $($or positive$)$  holomorphic sectional curvature? 
\item[$\bullet$] Can we find an analog of semi-positive $($or positive$)$ holomorphic sectional curvatures in algebraic geometry?
\end{itemize}
\end{prob}

\section{Pseudo-Effective Tangent Bundle} \label{Sec4}

Demailly-Peternell-Schenider in \cite{DPS94} established the structure theorem for compact K\"aher manifolds with nef tangent bundle. 
The nefness is generalized to the pseudo-effectivity in terms of singular hermitian metrics. 
The theory of singular hermitian metrics on vector bundles, which has been rapidly developed, gives a useful tool to study the pseudo-effectivity. 

This subsection is devoted to studying smooth projective varieties with pseudo-effective tangent bundle. 
We first recall our definition of pseudo-effective vector bundles and their characterizations. 
A vector bundle $E$ on a smooth projective variety $X$ is said to be \textit{positively curved} 
if $E$ admits a singular hermitian metric $h$ such that 
$\log |u|_{h^{\vee}}$ is a plurisubharmonic function for any local section $u$ of $E^\vee$, 
where $E^\vee$ is the dual vector bundle of $E$ and $h^\vee$ is the induced metric on $E^\vee$. 
See [HPS18, Definition 16.1] for the definition of singular hermitian metrics. 
Further $E$ is said to be \textit{pseudo-effective} 
if it satisfies one of the equivalent conditions in the following proposition:

\begin{prop}[{\cite[Proposition 3.1, Proposition 5.3]{BKK},  \cite[Subsection 2.3]{PT} }]
\label{chara}
Let $E$  be a vector bundle on a smooth projective variety $X$. 
Then the following conditions are equivalent: 
\begin{itemize}
\item[$(1)$] For any positive integer $m \in \mathbb{Z}_+$, there exists 
a singular hermitian metric $h_m$ on $\Sym ^m E$ such that 
$$
\deldel \log |u|^2_{h_m^{\vee}} \geq - \omega  
$$
for any local section $u$ of the $m$-th symmetric power $\Sym^m E^\vee$. 
Here $\omega$ is a fixed hermitian form on $X$ and $h_m^{\vee}$ is the induced metric on $\Sym^m E^\vee$. 
\item[$(2)$] There exists an ample line bundle $A$ such that $A \otimes \Sym^m{E}$ is generically globally generated for any integer $m > 0$ 
$($i.e., $A \otimes \Sym^m{E} $ is generated by global sections at a general point$)$.   
\item[$(3)$] Let $\mathcal{O}_{\mathbb{P}(E) }(1)$ be the hyperplane bundle on the projective space bundle $\mathbb{P}(E) \to X$. 
The non-nef locus  of  $\mathcal{O}_{\mathbb{P}(E) }(1)$ is not dominant over $X$. 
\end{itemize}
\end{prop}

Our definition requires that  the image of the non-nef locus of $\mathcal{O}_{\mathbb{P}(E) }(1)$ is properly contained in $X$, 
which is stronger than the pseudo-effectivity of $\mathcal{O}_{\mathbb{P}(E) }(1)$. 
Note that a smooth projective variety $X$ is isomorphic to the projective space 
if the tangent bundle $T_X$ is big in the following sense: 
the non-ample locus of  $\mathcal{O}_{\mathbb{P}(E) }(1)$ is not dominant over $X$ 
(see \cite[Corollary 6.7]{FM21}).  
See \cite{ELMNP06}, \cite{ELMNP09}, and \cite{Bou04} for non-nef loci and non-ample loci. 

The following theorem is a structure theorem of smooth projective varieties with pseudo-effective tangent bundle:


\begin{theo}[{\cite[Theorem 1.1]{HIM19}}]
\label{thm-psef}
Let $X$ be a smooth projective variety with pseudo-effective tangent bundle. 
Then $X$ admits a holomorphic MRC fibration $\phi: X \to Y$ 
to a smooth projective variety $Y$ with the following properties\,$:$
\begin{itemize}
\item[$\bullet$] The morphism $\phi: X \to Y$ is smooth. 
\item[$\bullet$] The image $Y$ admits a finite \'etale cover $A \to Y$ by an abelian variety $A$. 
\item[$\bullet$] A very general fiber $F$ of $\phi$ also has the pseudo-effective tangent bundle. 
\end{itemize}
Moreover, if $T_X$ is positively curved, 
then we obtain$:$
\begin{itemize}
\item[$\bullet$] The tangent bundle $T_X$ is decomposed into $T_X\cong T_{X/Y}  \oplus \phi^* T_Y$. 
\item[$\bullet$] The morphism $\phi: X \to Y$ is locally trivial. 
\end{itemize}
\end{theo}

Our structure theorem contains two essential differences from the case of nef tangent bundles in \cite{DPS94}. 
The first difference concerns  a splitting of the tangent bundle. 
In the case of $T_X$ being nef, 
the tangent bundle $T_X$ splits into $T_X \cong T_{X/Y} \oplus \phi^* T_Y$ and 
the subbundle $\phi^* T_Y \subset T_X$ is integrable. 
In particular, the splitting of $T_X$ implies that $\phi: X \to Y$ is locally trivial by \cite[Lemma 3.19]{Hor07}.
Under the weaker assumption of $T_X$ being pseudo-effective, 
the tangent bundle does not always split. 
Indeed, it follows from \cite[Proposition 4.2]{HIM19} that 
the tangent bundle $T_X$ of the projective space bundle $X:=\mathbb{P}(E) \to C$ over an elliptic curve $C$ 
is pseudo-effective, but does not split into the direct sum (in particular, it is not positively curved by the above theorem), 
where $E=\mathcal{O}_{C} \oplus \mathcal{O}_C(np)$, $p \in C$, and $n \in \mathbb{Z}_+$. 
Nevertheless, as of this moment, we have no counter-example to the local triviality of  $\phi: X \to Y$. 
We suggest the following problem: 

\begin{prob}\label{Fano}
Let $X$ be a smooth projective variety with pseudo-effective tangent bundle. 
Then is the MRC fibration $\phi : X \to Y$ locally trivial?
\end{prob}

The second difference concerns the positivity of fibers $F$ of $\phi:X \to Y$. 
It was proved in \cite{DPS94} that rationally connected manifolds with nef tangent bundle are always Fano manifolds. 
However, the same conclusion does not hold for pseudo-effective tangent bundles, 
since the tangent bundle of the Hirzebruch surfaces is pseudo-effective (see \cite[Proposition 4.5]{HIM19}). 
Then we suggest the following problem:

\begin{prob}[{\cite[Problem 3.13]{HIM19}}]
\label{Fano}
Let $X$ be a rationally connected manifold with pseudo-effective tangent bundle. 
Then is the anti-canonical divisor $-K_X$ big?
\end{prob}

Toward  structure theorems for compact K\"ahler manifolds $X$ with pseudo-effective tangent bundle, 
it was proved in  \cite[Theorem 3.12]{HIM19} 
that the Albanese map $X \to \Alb(X)$ satisfies a similar conclusion to Theorem \ref{thm-psef}. 
The next problem is to prove the projectivity of fibers. 

\begin{prob}\label{fiber}
Let $X$ be a compact K\"ahler manifold with pseudo-effective tangent bundle. 
After we replace $X$ with an appropriate finite \'etale covering of $X$, 
we consider the Albanese map $X \to \Alb(X)$. 
\begin{itemize}
\item[$\bullet$]  Can we prove that the anti-canonical divisor $-K_F$ is big or $F$  is projective?
\item[$\bullet$]  Is the Albanese map $X \to \Alb(X)$ an MRC fibration of $X$?
\end{itemize}
\end{prob}

It is not so easy to find examples of pseudo-effective tangent bundles. 
In \cite{HIM19}, we investigated surfaces with pseudo-effective tangent bundle 
by using the condition of Proposition \ref{chara} (2). 
In \cite{Par11}, examples of pseudo-effective tangent bundles are constructed 
in a similar way. 
On the other hand, H\"oring-Liu-Shao  in \cite{HLS20} also found examples of pseudo-effective tangent bundles 
by a different method based on VMRT (varieties of minimal rational tangents). 
Note that the definition of pseudo-effective vector bundles in \cite{HLS20} is weaker than our definition of Proposition \ref{chara}. 
The remaining problem in the classification of surfaces is as follows: 

\begin{prob}
$\bullet$ Does the blow-up of Hirzebruch surfaces at general four points have the pseudo-effective tangent bundle?
\\
$\bullet$ What can we say for the blow-up of Hirzebruch surfaces at special points?
\end{prob}

\section{Nef Anti-Canonical divisor} \label{Sec5}

This subsection is devoted to studying nef anti-canonical divisors.
The recent breakthrough in this direction 
is the structure theorem for smooth projective varieties with nef anti-canonical divisor by the works of \cite{Pau97, Zha96, Zha05, Cao19, CH19}, 
which can bee seen as an extension of the classical Beauville-Bogomolov decomposition to nef anti-canonical divisors. 
The structure theorem was proved even for compact K\"ahler manifolds by a generalized holonomy principle 
when the anti-canonical divisor admits a smooth hermitian metric with semi-positive curvature (see \cite{CDP15}). 
On the other hand, from the viewpoint of the minimal model program, 
the Beauville-Bogomolov decomposition was generalized to klt projective varieties 
by the works of \cite{GKP16, Dru18, GGK19, HP19}.

The most important remaining problem 
is to establish the structure theorem for klt pairs $(X, \Delta)$ with nef anti-canonical divisor $-(K_X+\Delta)$. 
Wang in \cite{Wang20} partially solved this problem when the regular locus $X_{\reg}$ of $X$ is simply connected and reduced the problem in the general case to some conjectures on fundamental groups of $X_{\reg}$. 
In the case $X$ being smooth, this problem was solved by \cite{CCM19} in the following form:

\begin{theo} [{\cite[Theorem 1.3]{CCM19}}]
\label{thm-nef}
Let $(X, \Delta)$ be a  klt pair with nef anti-canonical divisor $-(K_X +\Delta)$. 
Assume that $X$ is a smooth projective variety. 
Then there exists a holomorphic MRC  fibration $\phi: X \rightarrow Y$ with the following properties$:$
\begin{enumerate}
\item[$\bullet$] $Y$ is a smooth projective variety with numerically trivial canonical divisor. 
\item[$\bullet$] $\phi$ is locally trivial with respect to $(X, \Delta)$, i.e.,  for any small open set $U \subset Y$, we have the isomorphism
$$(\phi^{-1} (U), \Delta )\cong  U \times (X_y, \Delta_{X_y})$$
over $U \subset Y$.
Here $X_y$ is the typical fiber of $\phi$.
More strongly, there exits an isomorphism 
$$
(X_{\rm{univ}}, \Delta_{\rm{univ}}) \cong Y_{\rm{univ}} \times (X_{y}, \Delta_{X_y})
$$
over $Y_{\rm{univ}} $ and a representation $\rho : \pi_1(Y) \to {\rm{Aut}}(X_{y})$ 
such that $(X, \Delta)$ is isomorphic to the quotient of the right hand side by $\pi_1(Y) $. 
Here $X_{\rm{univ}}$ $($resp. $Y_{\rm{univ}}$$)$ is the universal covers of $X$ $($resp. $Y$$)$ and  
$\Delta_{\rm{univ}}$ is the pull-back of $\Delta $ to $X_{\rm{univ}}$. 
\end{enumerate} 

Moreover, together with the Beauville-Bogomolov decomposition, 
the universal cover $X_{\rm univ}$ of $X$ can be decomposed into 
the product of rationally connected manifolds, $\mathbb{C}^m$, Calabi-Yau manifolds, hyperk\"ahler manifolds. 
\end{theo}

The above structure theorem in the case of $\Delta=0$ was established
by the works of \cite{Pau97, Zha96, Zha05, Cao19, CH19} via the following steps: 
\begin{itemize}
\item[$(1)$] \cite{Pau97}: The fundamental group of $X$ is shown to have a polynomial growth 
by using the theory of Cheeger-Colding (see \cite{CC96}).

\item[$(2)$] \cite{Cao19}:  The problem is reduced to the case of $X$ being simply connected 
by studying of the Albanese map of $X$. 

\item[$(3)$] \cite{CH19}: The problem is solved for simply connected $X$
by studying MRC fibrations of $X$. 

\end{itemize}
It seems to be difficult to show that the fundamental group of $X$ have a polynomial growth directly for klt pairs in Theorem \ref{thm-nef}. 
Now we explain  a more direct approach suggested in \cite{CCM19}. 
The key point in the second and third steps is 
to show that the direct image sheaf $\phi_*(\hat A)$ satisfies a certain flatness,  
by applying the theory of positivity of direct images developed by \cite{Ber09, BP08, HPS18, PT}, 
where $\phi$ is the Albanese map (or an MRC fibration) of $X$ and  $\hat A$ is an appropriate relatively ample divisor. 
The direct image $\phi_*(\hat A)$ can be shown to be a trivial vector bundle when $X$ is simply connected. 
In \cite{CCM19}, without any assumptions on fundamental groups, 
we showed that $\phi_*(\hat A)$ admits a flat connection on a Zariski open set of $Y$ and it leads to the desired structure theorem.

It was revealed in \cite{CCM19} that almost all the arguments work for the nef relative anti-canonical divisor $-K_{X/Y}$ rather than the nef anti-canonical divisor $-K_X$. 
For example, when we consider a surjective morphism $\phi: X \to Y$ with nef relative anti-canonical divisor $-K_{X/Y}$, 
we can compare the Kodaira dimension and numerical Kodaira dimension of $-K_{X/Y}$ and $-K_{F}$ 
by applying positivity of direct images and the extension theorem in \cite{CDM17}. 
This study is motivated by Hacon-$\mathrm{M^{c}}$Kernan's question in \cite[Question 3.1]{HM07} (see \cite[Subsection 1.4]{EIM20} and \cite{EG19} for more details).  
Hence a structure theorem can also be  expected under the assumption for relative anti-canonical divisors. 

To explain that, let $\phi: X \dashrightarrow Y$ be an almost holomorphic map (not necessarily an MRC fibration) between smooth projective varieties, 
and let $\pi: \Gamma \to X$ be a resolution of the indeterminacy locus of $\phi: X \dashrightarrow Y$ 
with the corresponding morphism $\varphi: \Gamma \to Y$ in the following diagram: 
\begin{align*}
\xymatrix{
\Gamma \ar[rd]_\varphi  \ar[rr]^\pi & & X \ar@{-->}[ld]^\phi
\\
& Y. &}
\end{align*}
Assume that there exists an effective divisor $D$ on $\Gamma$ such that 
$(\Gamma, D)$ is a klt pair and $-\pi_*(K_{\Gamma/Y }+D)$ is a nef $\mathbb{Q}$-divisor.  
Note that this assumption is automatically satisfied 
when $(X, \Delta)$ is a klt pair with the nef anti-canonical divisor $-(K_X+\Delta)$ 
and $\phi: X \dashrightarrow Y$ is an MRC fibrations of $X$. 
Then, under this assumption, it can be shown that
a direct image sheaf of the form
$$
\varphi_*(-m(K_{\Gamma/Y} + D) + E +\hat A )
$$
satisfies a certain flatness  if $X$ is $\mathbb{Q}$-factorial (see \cite[Proposition 3.9]{CCM19} for the precise statement), 
where $E$ (resp.  $\hat A$) is an appropriate exceptional (resp. relatively ample) divisor.
This leads to the splitting $T_X=V_1 \oplus V_2$ such that 
$V_1$ coincides with $\pi_* T_{\Gamma/Y}$,  
if general fibers of $\phi: X \dashrightarrow Y$ are rationally connected (see \cite[Theorem 4.3]{CCM19}).   
Further, in this case, we can choose a holomorphic morphism $X \to Y'$
by replacing $Y$ with its another smooth birational model $Y'$. 
However, when $\phi: X \dashrightarrow Y$ is not an RC fibration, 
we do not know whether the same conclusion holds.
Then we have the following problem: 

\begin{prob}
Under the above assumption for the relative anti-canonical divisor, can we obtain a structure theorem?
Specifically, can we choose a holomorphic morphism $X \to Y'$ that is birationally equivalent to the original $\phi: X \dashrightarrow Y$?
\end{prob}

As explained in Section \ref{Sec4}, 
the structure theorem for nef tangent bundles in \cite{DPS94} is partially generalized to pseudo-effective tangent bundles. 
In the same spirit, it is natural and of interest to ask the following problem: 

\begin{prop}
What can we say for $($smooth$)$ projective varieties with pseudo-effective anti-canonical divisor?
\end{prop}

Toward the above problem, 
it seems to be the first step to consider a surjective morphism $\phi: X \to Y$ with the pseudo-effective relative anti-canonical $-K_{X/Y}$. 
We in  \cite{EIM20} systematically studied a relation between the geometric structure of $\phi$ and positivity conditions on $-K_{X/Y}$. 
As a result, we showed that the non-nef locus $\mathbb{B}_{-}(-K_{X/Y})$ is empty or dominant over $Y$ (see \cite{EIM20} for more details). 
This implies that the structure of $\phi: X \to Y$  is restricted even when $-K_{X/Y}$ is pseudo-effective.

The MRC fibration is defined only by $X$, 
and thus it does not reflect information on the boundary $\Delta$. 
On the other hand, 
the slope rationally connected quotient (sRC quotient for short), 
introduced in \cite{Cam16}, 
gives a generalization of MRC fibrations reflecting the boundary $\Delta$ (see \cite{Cam04, Cam16}).
The following conjecture concerns  a structure theorem for sRC quotients. 
This conjecture was solved when $(X, \Delta)$ is a log smooth surface in \cite[Theorem 1.6]{CCM19}, 
but it is still open in the general case.

\begin{conj}[{\cite[Conjecture 1.5]{CCM19}}]
\label{conj-sRC}
Let $(X, \Delta)$ be a klt pair such that $X$ is smooth and $- (K_X +\Delta)$ is nef. 
Then there exists an orbifold morphism $\rho : (X, \Delta)\rightarrow (R, \Delta_R)$ with the following properties$:$
\begin{enumerate}
\item[$\bullet$]  $(R, \Delta_R)$ is a klt pair such that $R$ is smooth and $c_1(K_R+ \Delta_R)=0$.

\item[$\bullet$]   General orbifold fibers $(X_r, \Delta_r )$ are slope rationally connected.

\item[$\bullet$]  The fibration is locally trivial with respect to pairs, namely, for any small open set $U \subset R$, we have the isomorphism over $U \subset Y$$:$ 
$$
(\rho^{-1} (U), \Delta)\cong  (U, \Delta_R |_U) \times (X_r, \Delta_{X_r}), 
$$
where $X_r$ is a general fiber of $\rho$.

\end{enumerate}
\end{conj}




\begin{thebibliography}{n}

\bibitem[AHZ18]{AHZ18}
A. Alvarez, G. Heier,  F. Zheng, 
\textit{On projectivized vector bundles and positive holomorphic sectional curvature}, 
Proc. Amer. Math. Soc. {\bf{146}} (2018), no. 7, 2877--2882.




 


\bibitem[BDPP13]{BDPP}
S. Boucksom, J.-P. Demailly, M. P\u{a}un, T. Peternell, 
\textit{The pseudo-effective cone of a compact K\"ahler manifold and varieties of negative Kodaira dimension}, 
J. Algebraic Geom. {\bf{22}} (2013), no. 2, 201--248. 

\bibitem[Bea83]{Bea83} 
 A. Beauville, 
\textit{Vari\'et\'es K\"ahleriennes dont la premi\`ere classe de Chern est nulle}, 
Geometry and analysis on complex manifolds, 39--50, 
J. Differential Geom. {\bf{18}} (1983), no. 4, 755--782 (1984). 

\bibitem[Ber09]{Ber09} 
B. Berndtsson,  
\textit{Curvature of vector bundles associated to holomorphic fibrations}, 
Ann. of Math. (2) {\bf{169}} (2009), no. 2, 531--560. 


\bibitem[BKK+15]{BKK}
T. Bauer, S. J Kov\'acs, A. K\"uronya, E. C. Mistretta, T. Szemberg, S. Urbinati, 
\textit{On positivity and base loci of vector bundles}, 
Eur. J. Math.  {\bf{1}}  (2015),  no. 2, 229--249.
 
 
\bibitem[Bou04]{Bou04}
S.  Boucksom, 
\textit{Divisorial Zariski decompositions on compact complex manifolds}, 
Ann. Sci. \'Ecole Norm. Sup. (4) {\bf{37}} (2004), no. 1, 45--76.

\bibitem[BP08]{BP08} 
B. Berndtsson, M. P\u{a}un, 
\textit{Bergman kernels and the pseudoeffectivity of relative canonical divisors}, 
Duke Math. J. {\bf{145}} (2008), no. 2, 341--378.
 
\bibitem[Cam92]{Cam92}
F. Campana, 
\textit{Connexit\'e rationnelle des vari\'et\'es de Fano}, 
Ann. Sci. \'Ecole Norm. Sup. (4) {\bf{25}} (1992), no. 5, 539--545. 

\bibitem[Cam04]{Cam04}
F. Campana, 
\textit{Orbifolds, special varieties and classification theory}, 
Ann. Inst. Fourier (Grenoble) 54 (2004), no. 3, 499--630. 

\bibitem[Cam16]{Cam16} 
F. Campana,
\textit{Orbifold slope rational connectedness}, 
available at arXiv:1607.07829v2. 

\bibitem[Cao19]{Cao19}
J. Cao, 
\textit{Albanese maps of projective manifolds with nef anticanonical divisors}, 
Ann. Sci. \'Ec. Norm. Sup\'er. (4),  {\bf 52} (2019), no. 5, 1137--1154. 


\bibitem[CC96]{CC96}
J. Cheeger, T. H. Colding, 
\textit{Lower bounds on Ricci curvature and the almost rigidity of warped products}, 
Ann. of Math. (2) {\bf{144}} (1996), no. 1, 189--237.

\bibitem[CCM21]{CCM19}
F. Campana, J. Cao, S. Matsumura, 
\textit{Projective klt pairs with nef anti-canonical divisor}, 
to appear in Algebr. Geom, 
available at arXiv:1910.06471v1. 

\bibitem[CDM17]{CDM17} 
J. Cao, J.-P. Demailly, S. Matsumura,  
\textit{A general extension theorem for cohomology classes on non reduced analytic subspaces}, 
Sci. China Math. {\bf{60}} (2017), no. 6, 949--962. 


\bibitem[CDP15]{CDP15} 
F. Campana, J.-P. Demailly, T. Peternell,  
\textit{Rationally connected manifolds and semipositivity of the Ricci curvature}, 
Recent advances in algebraic geometry, 71--91, 
London Math. Soc. Lecture Note Ser., {\bf{417}}, Cambridge Univ. Press, Cambridge, 2015. 


\bibitem[CH19]{CH19}
J. Cao, A. H\"oring, 
\textit{A decomposition theorem for projective manifolds with nef anticanonical divisor}, 
J. Algebraic Geom. {\bf 28} (2019), 567--597. 



\bibitem[CP91]{CP91}
F. Campana, T. Peternell, 
\textit{Projective manifolds whose tangent bundles are numerically effective}, 
Math. Ann., {\bf{289}} (1991), 169--187.





 
 
 \bibitem[DPS94]{DPS94}
J.-P. Demailly, T. Peternell, M. Schneider, 
\textit{Compact complex manifolds with numerically effective tangent bundles}, 
J. Algebraic Geom., {\bf{3}}, (1994), no.2, 295--345.

\bibitem[DPS01]{DPS01}
J-P. Demailly, T. Peternell, M. Schneider, 
\textit{Pseudo-effective line bundles on compact K\"ahler manifolds}, 
Internat. J. Math.  {\bf{12}}  (2001),  no. 6, 689--741.


\bibitem[Dru18]{Dru18}
S. Druel,
\textit{A decomposition theorem for singular spaces with trivial canonical class of dimension at most five}, 
Invent. Math. {\bf{211}} (2018), no. 1, 245--296.

 
 \bibitem[DT19]{DT19}
S. Diverio, S. Trapani, 
\textit{Quasi-negative holomorphic sectional curvature and positivity of the canonical divisor}, 
J. Differential Geom. {\bf{111}} (2019), no. 2, 303--314.
 
 
 \bibitem[EG19]{EG19}
S. Ejiri, Y. Gongyo, 
\textit{Nef anti-canonical divisors and rationally connected fibrations}, 
Compos. Math. {\bf{155}} (2019), no. 7, 1444--1456.
 
\bibitem[EIM20]{EIM20}
S. Ejiri, M. Iwai, S. Matsumura, 
\textit{On asymptotic base loci of relative anti-canonical divisors of algebraic fiber spaces}, 
available at arXiv:2005.04566v1. 

\bibitem[ELMNP06]{ELMNP06}
L. Ein, R. Lazarsfeld, M. Musta\c{t}\v{a}, M. Nakamaye, M. Popa,
\textit{Asymptotic invariants of base loci}, 
Ann. Inst. Fourier (Grenoble) {\bf 56} (2006), no. 6, 1701--1734. 

\bibitem[ELMNP09]{ELMNP09} 
L. Ein, R. Lazarsfeld, M. Musta\c{t}\v{a}, M. Nakamaye, M. Popa, 
\textit {Restricted volumes and base loci of linear series},
Amer. J. Math. {\bf 131} (2009), no.3, 607--651.




\bibitem[FM21]{FM21}
M. Fulger, T. Murayama
\textit{Seshadri constants for vector bundles}, 
to appear in Journal of Pure and Applied Algebra, {\bf{225}}, Issue 4, (2021). 


\bibitem[GGK19]{GGK19}
D. Greb, H. Guenancia,  S. Kebekus, 
\textit{Klt varieties with trivial canonical class: holonomy, differential forms, and fundamental groups}, 
Geom. Topol. {\bf{23}} (2019), no. 4, 2051--2124.

 
\bibitem[GHS03]{GHS03}
T. Graber, J. Harris, J. Starr, 
\textit{Families of rationally connected varieties}, 
J. Amer. Math. Soc. {\bf{16}} (2003), no. 1, 57--67.


\bibitem[GKP16]{GKP16}
D. Greb, S. Kebekus, T. Peternell, 
\textit{Singular spaces with trivial canonical class}, 
Minimal models and extremal rays (Kyoto, 2011), 67--113, Adv. Stud. Pure Math., {\bf{70}}, Math. Soc. Japan, [Tokyo], 2016. 








\bibitem[HC20]{HC20}
A. Chaturvedi, G. Heier, 
\textit{Hermitian metrics of positive holomorphic sectional curvature on fibrations}, 
Math. Z. {\bf{295}} (2020), no. 1-2, 349--364.

\bibitem[HIM21]{HIM19}
G. Hosono, M. Iwai, S. Matsumura, 
\textit{On projective manifolds with pseudo-effective tangent bundle},  
to appear in J. Inst. Math. Jussieu, 
available at DOI: https://doi.org/10.1017/S1474748020000754

\bibitem[Hit75]{Hit75}
N. Hitchin, 
\textit{On the curvature of rational surfaces}, 
Differential geometry (Proc. Sympos. Pure Math., Vol. XXVII, Part 2, Stanford Univ., Stanford, Calif., 1973), 
pp. 65--80. Amer. Math. Soc., Providence, R. I., 1975.



\bibitem[HLS20]{HLS20}
A. H\"oring, J. Liu, F. Shao, 
\textit{Examples of Fano manifolds with non-pseudoeffective tangent bundle}, 
available at arXiv:2003.09476v1. 


\bibitem[HLWZ18]{HLWZ18}
G. Heier, S. S. Y. Lu, B. Wong, F. Zheng, 
\textit{Reduction of manifolds with semi-negative holomorphic sectional curvature}, 
Math. Ann. {\bf{372}} (2018), no. 3-4, 951--962.

\bibitem[HM07]{HM07}
C.D. Hacon,  J. $\mathrm{M^{c}}$Kernan, 
\textit{On Shokurov's rational connectedness conjecture}, 
Duke Math. J. {\bf{138}} (2007), no. 1, 119--136. 

\bibitem[H\"or07]{Hor07}
A. H\"oring, 
\textit{Uniruled varieties with split tangent bundle}, 
Math. Z., {\bf{256}} (2007), no.3, 465--479.




\bibitem[HP19]{HP19}
A. H\"oring, P. Peternell, 
\textit{Algebraic integrability of foliations with numerically trivial canonical divisor}, 
Invent. Math. 216 (2019), no. 2, 395--419.


\bibitem[HPS18]{HPS18}
C. Hacon, M. Popa, C. Schnell, 
\textit{Algebraic fiber spaces over abelian varieties: around a recent theorem by Cao and P{\v{a}}un}, 
Local and global methods in algebraic geometry, 143--195, 
Contemp. Math., {\bf{712}}, Amer. Math. Soc., Providence, RI, 2018.


\bibitem[HSW81]{HSW81}
A. Howard, B. Smyth, H. Wu, 
\textit{On compact K\"ahler manifolds of nonnegative bisectional curvature I and II}, 
Acta Math. {\bf{147}} (1981), no. 1-2, 51--70. 


\bibitem[HW20]{HW15}
G. Heier, B. Wong, 
\textit{On projective K\"ahler manifolds of partially positive curvature and rational connectedness}, 
Doc. Math. {\bf{25}} (2020), 219--238






\bibitem[Iwa18]{Iwa}
M. Iwai, 
\textit{Characterization of pseudo-effective vector bundles by singular hermitian metrics}, 
to appear in Michigan Math. J., 
available at arXiv:1804.02146v2. 


\bibitem[KoMM92]{KoMM92}
J. Koll\'ar, Y. Miyaoka, S. Mori, 
\textit{Rationally connected varieties}, 
J. Algebraic Geom. {\bf{1}} (1992), no. 3, 429--448.


\bibitem[Mat13]{Mat13}
S. Matsumura, 
\textit{Asymptotic cohomology vanishing 
and a converse to the Andreotti-Grauert theorem on surfaces}, 
Ann. Inst. Fourier (Grenoble) {\bf{63}} (2013), no. 6, 2199--2221.


\bibitem[Mat20]{Mat18a}
S. Matsumura, 
\textit{On the image of MRC fibrations of projective manifolds with semi-positive holomorphic sectional curvature}, 
Pure Appl. Math. Q.  {\bf{16}}, No. 5 (2020), pp. 1443--1463.
published online at https://dx.doi.org/10.4310/PAMQ.2020.v16.n5.a4. 


\bibitem[Mat21]{Mat18b}
S. Matsumura, 
\textit{On projective manifolds with semi-positive holomorphic sectional curvature}, 
to appear in Amer. J. Math., 
available at arXiv:1811.04182v1. 

\bibitem[Mok88]{Mok88}
N. Mok,
\textit{The uniformization theorem for compact K\"ahler manifolds of nonnegative holomorphic bisectional curvature}, 
J. Differential Geom. {\bf{27}} (1988), no. 2, 179--214.



\bibitem[Mor79]{Mor79}
S. Mori, 
\textit{Projective manifolds with ample tangent bundles}, 
Ann. of Math. (2) {\bf{110}} (1979), no. 3, 593--606. 


\bibitem[Nak04]{Nak04}
N. Noboru, 
\textit{Zariski-decomposition and abundance}, 
MSJ Memoirs, {\bf{14}}. Mathematical Society of Japan, Tokyo, 2004. 
xiv+277 pp. ISBN: 4-931469-31-0. 


\bibitem[Par11]{Par11}
M. Paris, 
\textit{Quelques aspects de la positivit\'e du fibr\'e tangent des vari\'et\'es projectives complexes}, 
available at  https://tel.archives-ouvertes.fr/tel-00552308. 



\bibitem[P\u{a}u97]{Pau97} 
M. P\u{a}un, 
\textit{Sur le groupe fondamental des vari\'et\'es k\"ahl\'eriennes 
compactes \`a classe de Ricci num\'eriquement effective}, 
C. R. Acad. Sci. Paris S\'er. I Math. {\bf{324}} (1997), no. 11, 1249--1254.




 
 \bibitem[PT18]{PT} 
M. P\u{a}un, S. Takayama, 
\textit{Positivity of twisted relative pluricanonical divisors and their direct images}, 
J. Algebraic Geom. {\bf{27}} (2018), 211--272.

\bibitem[Rau15]{Rau15}
H. Raufi,
\textit{Singular hermitian metrics on holomorphic vector bundles},
Ark. Mat. {\bf{53}} (2015), no. 2, 359--382.


\bibitem[SY80]{SY80}
Y.-T. Siu, S.-T. Yau, 
\textit{Compact K\"ahler manifolds of positive bisectional curvature}, 
Invent. Math. {\bf{59}} (1980), no. 2, 189--204.


\bibitem[TY17]{TY17}
V. Tosatti, X. Yang, 
\textit{An extension of a theorem of Wu-Yau}, 
J. Differential Geom. {\bf{107}} (2017), no. 3, 573--579.
 
 
 \bibitem[Wan20]{Wang20}
J. Wang, 
\textit{Structure of projective varieties with nef anticanonical divisor: the case of log terminal singularities}, 
available at arXiv:2005.05782v2. 

\bibitem[Wu20]{Wu20}
X. Wu, 
\textit{Pseudo-effective and numerically flat reflexive sheaves}, 
available at arXiv:2004.14676v2. 


\bibitem[WY16]{WY16}
D. Wu, S.-Y. Yau, 
\textit{Negative holomorphic curvature and positive canonical divisor}, 
Invent. Math. {\bf{204 (2016)}}, no. 2, 595--604.




\bibitem[Yan18a]{Yan18a}
X. Yang, 
\textit{RC-positivity, rational connectedness and Yau's conjecture}, 
Camb. J. Math.  {\bf{6}} (2018), 183--212. 




\bibitem[Yan18b]{Yan18b}
X. Yang, 
\textit{RC-positive metrics on rationally connected manifolds}, 
Forum Math. Sigma {\bf{8}} (2020), Paper No. e53, 19 pp.

\bibitem[Yan19]{Yan19}
X. Yang, 
\textit{A partial converse to the Andreotti-Grauert theorem}, 
Compos. Math. {\bf{155}} (2019), no. 1, 89--99. 
 

\bibitem[Yau82]{Yau82}
S.-T. Yau, 
\textit{Problem section}, 
Seminar on Differential Geometry, 669--706, 
Ann. of Math. Stud., {\bf{102}}, Princeton Univ. Press, Princeton, N.J, (1982).

\bibitem[Zha96]{Zha96} 
Q. Zhang, 
\textit{On projective manifolds with nef anticanonical divisors},  
J. Reine Angew. Math. {\bf{478}} (1996), 57--60. 

\bibitem[Zha05]{Zha05} 
Q. Zhang, 
\textit{On projective varieties with nef anticanonical divisors},  
Math. Ann. {\bf{332}} (2005), no. 3, 697--703.



\end{thebibliography}
\end{document}